\documentclass[11pt]{mypaper}
\usepackage[all]{xy}
\xyoption{poly}
\acinq
\aquatre
\setlist{nosep}

\newcommand{\F}{\mathcal{F}}
\newcommand{\C}{\mathcal{C}}

\renewcommand{\P}{\mathcal{P}}
\newcommand{\M}{\mathcal M}
\newcommand{\K}{\mathcal K}

\newcommand{\abar}{\overline a}

\newcommand{\bbar}{\overline b}

\newcommand{\ybar}{\overline y}
\newcommand{\Vbar}{\overline V}

\newcommand{\I}{\mathcal I}
\newcommand{\D}{\mathcal D}
\newcommand{\E}{\mathcal E}
\renewcommand{\SS}{\mathcal S}

\newcommand{\ve}{\varepsilon}
\newcommand{\ie}{{\normalfont ie}}

\newcommand{\pr}{\mathbb{P}}

\newcommand{\un}{\mathbf{1}}

\newsavebox{\ddtbox}
\savebox{\ddtbox}{$\frac{\mathrm{d}}{\mathrm{d}t}$}
\newcommand{\ddt}{\usebox{\ddtbox}}

\newcommand{\iset}{independence set}
\newcommand{\dset}{dependence set}
\newcommand{\config}{configuration}
\newcommand{\wconfig}{weighted \config}
\newcommand{\one}{{\normalfont I}}
\newcommand{\two}{{\normalfont II}}

\begin{document}

\begin{center}
  \huge
  \textbf{Configured spaces and their Möbius polynomials}\par\medskip
  \Large Samy Abbes\\
  \large Université Paris Cité --- IRIF CNRS UMR 8243\\
  \normalsize\texttt{abbes@irif.fr}\\
  August 2025
\end{center}

\begin{abstract}
  An arbitrary dependence structure between a finite family of events of a probability space defines a hypergraph structure. We study the converse operation, starting from a hypergraph structure, to determine a canonical probability space with events indexed by the vertices and that satisfy the prescribed dependence structure; we also require that events that are not mutually independent must be mutually exclusive. Such a hypergraph structure, we call a configuration.

  Requiring furthermore that all events are given the same positive probability~$t$, converts this problem to the study of a family of Möbius polynomials, which we call the relative Möbius polynomials (or independence polynomial) of the configuration. This study was very much already complete in the case where the configuration actually derives from a graph structure, but remained open in the general case.

  We show the existence and uniqueness of an optimal real~$t$, optimal in the sense that the events cover a subset of maximal probability. We base our study on a new formula for the derivative of the Möbius polynomial.

\end{abstract}

\section{Introduction}
\label{sec:introduction}

Let $(\Omega,\F,\pr)$ be a probability space, and let $(a_1,\dots,a_n)$ be a family of events, all of them of positive probability, and with the following property: for every subset $A\subseteq [n]=\{1,\dots,n\}$, either the family $(a_i)_{i\in A}$ is exclusive modulo zero, ie, $\pr\bigl(\bigcap_{i\in A}a_i\bigr)=0$, or it is an independent family of events, ie, $\pr(\bigcap_{i\in B}a_i)=\prod_{i\in B}\pr(a_i)$ for every subset $B\subseteq A$. This whole data makes what we call a \emph{configured space}.

This configured space induces a hypergraph structure on $V=[n]$, namely it defines the two families $\I=\{x\in\P(V)\tq \text{$x$ is independent}\}$ and $\D=\{x\in\P(V)\tq \text{$x$ is mutually exclusive}\}$. Let $\leq$ denote the inclusion order on~$\P(V)$. Obviously, $\I$~is $\leq$-downward closed, whereas $\D$ is $\leq$-upward closed, and $\I$ contains all singletons.

Conversely, consider a hypergraph structure $(V,\I)$ where $V$ is a finite set of vertices, and $\I$ is a family of subsets of~$V$, such that $\I$ contains all singletons $\{a\}$ for $a\in V$ and is downward closed with respect to inclusion in~$\P(V)$. Put also $\D=\P(V)\setminus\I$. The data $(V,\I,\D)$ forms what we call a \emph{configuration}. 

A first question is the existence of a probability space $(\Omega,\F,\pr)$ and of a mapping $\varphi:V\to\F$ with the following properties:
\begin{inparaenum}[1)]
  \item $\pr(\varphi(a))>0$ for all $a\in V$; and
  \item  $\varphi(x)$~is an independent family of events if $x\in\I$; and
  \item  $\varphi(x)$~is mutually exclusive modulo zero if $x\in\D$.
  \end{inparaenum}

  A second question consists in replacing the first item above, ie, $\pr(\varphi(a))>0$ for all $a\in V$, by the requirement that $\pr(\varphi(a))=t$ for all $a\in V$ and for some positive real~$t$. Finally, a third question consists in adding the following requirement:  the whole collection of events $(\varphi(a))_{a\in V}$ forms a covering modulo zero of~$\Omega$, \ie, $\pr(\bigcup_{a\in V}\varphi(a))=1$.

Although this topic does not seem to have been systematically explored so far, it turns out that the above questions have been given an answer for some particular cases, namely when the abstract configuration $(V,\I,\D)$ is induced by a graph $(V,D)$, in the sense that $\I$ is the set of independence cliques of the graph. In~\cite{shearer85}, Shearer uses this construction to establish optimal bounds for the Lovász Local Lemma; many works have been inspired from Shearer's paper, the main theme being to establish optimal regions where the independence polynomials of graphs with bounded degree do not vanish. Extensions of results in the same vein have been carried on to hypergraphs~\cite{galvin23}. The present work departs from this vein of research as we are interested in studying precise independence polynomials, and not in obtaining general results valid for hypergraphs of a given maximal degree. Furthermore, our probabilistic motivations narrows our study of the independence polynomials on the negative part of the line only, and more precisely on the $(-1,0)$ interval only. Actually, the Möbius polynomials $\mu(t)$ that we study are related to the independence polynomials $\theta(t)$ by $\mu(t)=\theta(-t)$, and our study related to the Möbius polynomials is thus located on the interval $(0,1)$.

Here is how the Möbius polynomial appears, and this will also introduce the major difference between the general case and the case of graphs. Assume given a configured space $(\Omega,\F,\pr)$ with configuration $(V,\I,\D)$ and such that $\pr(a)=t>0$ for all $a\in V$; for simplicity, we drop the mapping~$\varphi$, and simply identify vertices with the events of $\Omega$ that they represent. Let the \emph{rest} be defined by $R(t)=1-\pr\bigl(\bigcup_{a\in V}a)$, so that the configuration forms a covering modulo zero if and only if $R(t)=0$. By the inclusion-exclusion principle, we have:
\begin{gather}
  \label{MMeq:1}
  R(t)=1-\sum_{x\in\P(V),\ x\neq\emptyset}(-1)^{|x|+1}\pr\Bigl(\bigcap x\Bigr)
\end{gather}

Now every non empty subset $x\in\P(V)$ is either a family of independent events, if $x\in\I$, in which case $\pr\bigl(\bigcap x)=\prod_{a\in x}\pr(a)=t^{|x|}$; or a family of mutually exclusive events modulo zero, if $x\in\D$, in which case $\pr(\bigcap x)=0$. Hence:
\begin{gather}
  \label{MMeq:2}
  R(t)=\sum_{x\in\I}(-1)^{|x|}t^{|x|}
\end{gather}

Therefore, a necessary condition for the existence of $t$ satisfying our requirements is that $R(t)\geq0$, where $R(t)$ is now the polynomial defined by~\eqref{MMeq:2}, and which only depends on the abstract configuration. This polynomial is the Möbius polynomial of the configuration. $R(0)=1$ hence $R>0$ in a neighborhood of~$0$, but we shall see that the condition $R(t)\geq0$ is not sufficient for $t$ to define a valid probability~$\pr_t$, ie, satisfying $\pr_t(a)=t$ and respecting the independence and exclusion conditions of the abstract configuration. Nevertheless, we shall see that there is indeed an interval $[0,\delta]$ for which all $t\in(0,\delta]$ define a valid configured space with $\pr_t(a)=t$ for all $a\in V$.

A secondary question that arises is the monotony of $R(t)$ on $[0,\delta]$. Of course, as a polynomial of arbitrary degree, it has no reason to be monotonic on the interval $(0,1)$. From the probabilistic point of view, $1-R(t)$ is the probability of the union $\bigcup V$. When $t$ increases, every individual probability $\pr_t(a)=t$ increases, but this is not sufficient in general to guarantee that $t\mapsto\pr_t(\bigcup V)$ increases, and thus that $t\mapsto R(t)$ decreases. Nevertheless, we shall prove this fact indeed, and for this we give an expression for the derivative of $R(t)$ which does not seem to have been noted so far, and which is of interest on its own.

As a result, if the configuration can define a covering of its configured space for some real~$t>0$, then it must be the case that $t=\delta$ and that $R(\delta)=0$.

Now if we examine the case of a configuration induced by a graph, several questions already have an answer. The major advantage in this situation is that the polynomial $R(t)$ has a combinatorial interpretation, as follows. Let $\M$ be the trace monoid generated by $(V,D)$. That is, introducing the binary relation $I=(V\times V)\setminus D$, $\M$~is the monoid defined by generators and relations by: $\M=\langle V\ |\ \forall (a,b)\in I\quad ab=ba\rangle$, and let $G(t)=\sum_{x\in\M}t^{|x|}$ be the generating series of~$\M$. Then a result of Cartier and Foata in~\cite{cartier69} is the following identity: $R(t)G(t)=1$, hence $R(t)=\frac1{G(t)}$. But~$G(t)$ is a power series with non negative coefficients. Therefore its radius of convergence is a pole of the series, hence a zero of the polynomial $R(t)$, and it can be proved that this real is indeed the real~$\delta$. Therefore the covering is always possible, realized with~$\delta$, and it is also immediate that $t\mapsto R(t)$ is decreasing on $[0,\delta]$, again by using the fact that all the coefficients of $G(t)$ are non negative.

Nothing could be better if this approach could be extended to the case of a general abstract configuration. But, this is impossible, because simple examples show that, in general, the inverse of the polynomial $R(t)$ is not a series with non negative coefficients. Henceforth, a nice combinatorial interpretation of this kind for $\frac1{R(t)}$ is impossible in general. For example, take $V=\{1,2,3\}$ and $\I$ as all subsets of~$V$, but $V$~itself. Then $R(t)=1-3t+3t^2$ has no real root, hence its inverse cannot have only non negative coefficients. 

\paragraph{Contributions.}

In this paper we develop some elements for a theory of configurations and of configured spaces. In particular, beside the Möbius polynomial of a configuration, which is a long-studied object in hypergraph theory, we introduce relative configurations and relative Möbius polynomials, and we show that the configuration is best studied with the help of this whole family of polynomials, rather than with its single Möbius polynomial. In particular, we give an expression for the derivative of the Möbius polynomial in terms of the relative polynomials. This allows us to determine the simple form $(0,\delta]$ of the reals making a configuration eligible for a representation as a configured space. We therefore identify two types of configurations, those of type \one\ for which the covering is possible, and those of type \two\ for which the covering is not possible. As we already noticed, configurations induced by graphs are all of type~\one. We introduce the example of the star configurations, a family of examples containing configurations of both types.

Finally we gather in the last section some questions that we think are of interest for further study.

\section{Free models of probability spaces through events}
\label{sec:free-models-prob}

This preliminary section gathers some folklore elementary results. Let $a_1,\ldots,a_n$ be $n$ distinct symbols, thought of as events of some probability space to be constructed. For each subset $A\subseteq [n]=\{1,\dots,n\}$, let $p_A$ be a non negative real, thought of as the probability of $\bigcap A=\bigcap_{i\in A}a_i$. What are necessary and sufficient conditions on the family $(p_A)_{A\in\P([n])}$ for the existence of a probability space $(\Omega,\F,\pr)$ and a mapping $\varphi:[n]\to\F$ such that $\pr\bigl(\bigcap\varphi( A)\bigr)=p_A$ for every $A\in\P([n])$~?

We will give below a necessary and sufficient condition, together with the construction of canonical probability space $(\Omega,\F,\pr)$ and of a canonical mapping $\varphi:[n]\to\F$. This construction is of the type frequently encountered in the probabilistic method literature, see for instance~\cite{alon00}.

We start with an easier question. Let $V=\{a_1,\ldots,a_n\}$. For each $a\in V$, introduce a distinct symbol~$\abar$, let $\Vbar=\{\abar_1,\dots,\abar_n\}$ and let $\K$ be the set of words formed with symbols in $V\cup\Vbar$, and such that for each $i\in[n]$, there is either an occurrence of~$a_i$, or an occurrence of~$\abar_i$, or none. Furthermore, words are formed regardless of the order of occurrences of symbols. For instance with $V=\{a,b\}$, and grouping words by increasing size, this yields:
\begin{gather*}
  \K=\{\ve,\quad a,b,\abar,\bbar,\quad ab,a\bbar,\abar b,\abar \bbar\}
\end{gather*}
with $\ve$ the empty word. 

We equip $\K$ with the prefix ordering~$\leq$. Recall that the order of letters is irrelevant, hence two words $x$ and $y$ satisfy $x\leq y$ if and only every symbol of $x$ occurs in~$y$. Let $\E$ be the set of words of maximal size, which is $n=|V|$. The set $\E$ is also the set of maximal elements in~$(\K,\leq)$.

Let $(\Omega,\F,\pr)$ be a probability space and let $\varphi:V\to\F$ be a mapping. We extend $\varphi$ to a mapping still denoted by $\varphi:\K\to\F$ by setting for each $x\in\K$, say $x=x_1\dots x_k$ with $x_i\in V\cup\Vbar$:
\begin{gather}
  \label{MMeq:4}
  \varphi(x)=\bigcap_{1\leq i\leq k}X_i
\end{gather}
where $X_i=\varphi(a)$ if $x_i=a\in V$ and $X_i=\overline{\varphi(a)}$, the complementary of~$\varphi(a)$ in~$\Omega$, if $x_i=\abar\in\Vbar$. We also put $\varphi(\ve)=\Omega$.

Assume given a non negative real $p_x$ for every $x\in\E$. We say that $\varphi$ realizes the family $(p_x)_{x\in\E}$ if:
\begin{gather}
  \label{MMeq:5}
  (x\in\E)\quad \pr\bigl(\varphi(x)\bigr)=p_x
\end{gather}

\begin{proposition}
  \label{MMprop:1}
  There exists a probability space $(\Omega,\F,\pr)$ and a realization $\psi:V\to\F$ if and only if $\sum_{x\in\E}p_x=1$.

  If $\sum_{x\in\E}p_x=1$, then a canonical realization is the finite probability space $(X,m)$ where $X=\E$ equipped with the probability distribution $m=(p_x)_{x\in\E}$, and with the realization $\varphi:V\to\P(X)$, $a\mapsto\varphi(a)=\{x\in\E\tq a\leq x\}$. Furthermore, for every $x\in\K$, the event $\varphi(x)$ has the probability:
  \begin{gather}
    \label{MMeq:3}
    \pr\bigl(\varphi(x)\bigr)=\sum_{y\in\E\tq x\leq y}p_y
  \end{gather}

If $(\Omega,\F,\pr)$ provides another realization, then there is a canonical random variable $U:\Omega\to X$ with law~$m$.
\end{proposition}

\begin{proof}
  If a realization exists, then the family of events $(\psi(x))_{x\in\E}$ is a partition of~$\Omega$, whence the necessity of $\sum_{x\in\E}p_x=1$.

  Conversely, if $\sum_{x\in\E}p_x=1$, then $(X,m)$ defined as in the statement is a finite probability space by construction. Let $\varphi:V\to\F$ be defined as in the statement, and let $\varphi:\K\to\F$ be its extension defined as in~\eqref{MMeq:4}. Then, for $x\in\K$,  say $x=a_1\dots a_k$ with $a_i\in V\cup \Vbar$, one has $\varphi(x)=\{y\in\E\tq \forall i=1,\dots,k\quad a_i\leq y\}=\{y\in\E\tq x\leq y\}$. It implies on the one hand that $\varphi(x)=\{x\}$ for $x\in\E$, hence that $\varphi$ is a realization; and on the other hand that~\eqref{MMeq:3} holds for every $x\in\K$.

  Finally, if $(\Omega,\F,\pr)$ and $\psi$ provide another realization, the canonical random variable $U:\Omega\to\E$ is given by $U(\omega)=x_1\dots x_n$ where $x_i=a_i$ if $\omega\in\psi(a_i)$ and $x_i=\abar_i$ if $\omega\notin\psi(a_i)$. The law of $U$ is indeed given, for $x\in\E$, by $\pr(U=x)=\pr\bigl(\psi(x)\bigr)=p_x=m(x)$.
\end{proof}

We now come back to the original question of this section. Let $\K^+$ be the set of positive words of~$\K$, ie, the words with only occurrences of symbols in~$V$, and no occurrences of symbols in~$\Vbar$, and let $(q_x)_{x\in\K^+}$ be a family of non negative reals with $q_\ve=1$. For a probability space $(\Omega,\F,\pr)$ and a mapping $\varphi:V\to\F$, we say that $\varphi$ is a realization of $(q_x)_{x\in\K^+}$ if $\pr\bigl(\varphi(x)\bigr)=q_x$ for every $x\in\K^+$, where the event $\varphi(x)$ is defined as above in~\eqref{MMeq:4} by $\varphi(x)=\bigcap_{1\leq i\leq k}\varphi(x_i)$, for $x=x_1\dots x_k$ with $x_i\in V$.

Assume that such a realization exists.
By the inclusion-exclusion formula:
\begin{gather}
  \label{MMeq:7}
  \pr\bigl(\varphi(\abar_1\dots\abar_n)\bigr)=\sum_{x\in\K^+}(-1)^{|x|}q_x
\end{gather}
And more generally, for every $x\in\K^+$, if $\ybar$ is the only word which contains only occurrences of elements of~$\Vbar$ and such that $x\ybar\in\E$, one has:
\begin{gather}
  \label{MMeq:8}
  \pr\bigl(\varphi(x\ybar)\bigr)=\sum_{z\in\K^+\tq x\leq z}(-1)^{|z|-|x|}q_x
\end{gather}
The formula~\eqref{MMeq:7} is a particular case of~\eqref{MMeq:8}, obtained with $x=\ve$ in~\eqref{MMeq:8}.

Hence for instance with $n=3$, dropping the notation~$\varphi$ and identifying the events with their probabilities for a more concise notation:
\begin{align*}
  \abar_1\abar_2\abar_3&=1-a_1-a_2-a_3+a_1a_2+a_1a_2+a_2a_3-a_1a_2a_3\\
  a_1\abar_2\abar_3&=a_1-a_1a_2-a_1a_3+a_1a_2a_3\\
  a_1a_2\abar_3&=a_1a_2-a_1a_2a_3
\end{align*}

We conclude that, if the realization of $(q_x)_{x\in\K^+}$ exists, then the values of $p_x=\pr\bigl(\varphi(x)\bigr)$ for $x\in\E$ are all determined, independently of the space~$(\Omega,\F,\pr)$, and are given by:
\begin{gather}
  \label{MMeq:9}
(x\in\E)\quad  p_x=\sum_{z\in\K^+\tq y\leq z}(-1)^{|z|-|y|}q_z
\end{gather}
where $y$ is the subword of $x$ that contains all its positive symbols, ie, belonging to~$V$.

The sum of the family $(p_x)_{x\in\E}$ thus defined is:
\begin{align*}
\sum_{(y,z)\in\K^+\times\K^+\tq y\leq z}(-1)^{|z|-|y|}q_z=\sum_{z\in\K^+}(-1)^{|z|}q_z\Bigl(\sum_{y\in\K^+}(-1)^{|y|}\un(y\leq z)\Bigr)=1
\end{align*}
Hence, in view of Proposition~\ref{MMprop:1}, we arrive at the following statement; note that the property~\eqref{MMeq:6} below is independent of the realization, and only depends on the family~$(q_x)_{x\in\K^+}$.

\begin{proposition}
  \label{MMprop:2}
  For any family $(q_x)_{x\in\K^+}$ with $q_\ve=1$, there exists a realization of $(q_x)_{x\in\K^+}$ if and only if the reals $p_x$ defined by~\eqref{MMeq:9} for $x\in\E$ are all non negative.

  If this is the case, a canonical realization is given by the finite probability space $(X,m)$ with $x=\E$ and $m(x)=p_x$ for $x\in\E$, with the mapping $\varphi:V\to\P(X)$ given by $\varphi(a)=\{x\in\E\tq a\leq x\}$.

  If $(\Omega,\F,\pr)$ provides another realization of $(q_x)_{x\in\K^+}$, then there is a canonical random variable $U:\Omega\to X$ with law~$m$. The collection $(\varphi(a))_{a\in V}$ is a covering modulo zero of\/ $\Omega$ if and only if:
  \begin{gather}
    \label{MMeq:6}
    \sum_{x\in\K^+}(-1)^{|x|}q_x=0
  \end{gather}
\end{proposition}

\section{Configurations}
\label{sec:configurations}

\subsection{Sets}
\label{sec:sets}

\begin{definition}[configurations, nubs]
  \label{def:2}
  A \emph{\config} is a triple $(V,\I,\D)$ where $V$ is a finite set of \emph{vertices}, and $\I$ and $\D$ are two subsets of $\P(V)$ such that:
  \begin{inparaenum}[1)]
  \item $\I$~is non empty and downward closed in~$(\P(V),\subseteq)$;
  \item $\I$~contains all singletons $\{a\}$ for $a\in V$; and
  \item $\D=\P(V)\setminus\I$. 
  \end{inparaenum}
  
  Elements of $\I$ are called \iset s and elements of $\D$ are called \dset s. 
\end{definition}

\begin{remark}
  \label{rem:1}
  \begin{enumerate}
  \item Of course, the \config\  $(V,\I,\D)$ can be specified by $V$ and $\D$ only, setting $\I=\P(V)\setminus\D$, provided that $\D$ is upward closed  in $\P(V)$ and contains no singleton $\{a\}$ for $a\in V$.
  \item The possibility that $\D=\emptyset$ is not ruled out.
  \item The \emph{trivial} \config\ is $(\emptyset,\{\emptyset\},\emptyset)$, which we simply denote by~$\emptyset$.
  \end{enumerate}
\end{remark}

\textbf{Notations~:} Instead of the usual symbol~$\subseteq$, we shall use the symbol $\leq$ for denoting the inclusion ordering between subsets of~$V$, hence making $(\I,\leq)$ and $(\D,\leq)$ two partially ordered sets. We shall denote the empty set by~$\ve$, the singleton $\{a\}$ for $a\in V$ simply by~$a$, and any element of $\P(V)$ with elements $a_1,\ldots,a_k$ by the word $a_1\dots a_k$, where it is understood that the order of letters is irrelevant.

\begin{definition}
  \label{def:6}
The \emph{nubs} of a configuration $\C=(V,\I,\D)$ are the minimal elements of $(\D,\leq)$.
\end{definition}

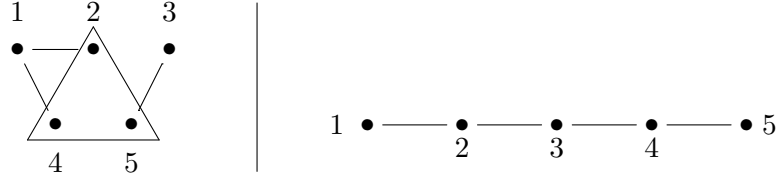
\begin{figure}
$$
\begin{array}{c|c}
  \begin{xy}
    0;/r1cm/:
    (.5,0)="A",*+{\bullet}
    ,{\PATH ~={**\dir{-}}
      '(0,1)="C"*+{\bullet} (1,1)="D"*+{\bullet}}
,0
,{(1.5,0)="B"*+{\bullet}    \ar@{-} (2,1)="E"*+{\bullet}}
,(1,.3),{\xypolygon3{}}
,"C"+(0,.5)*{1}
,"D"+(0,.5)*{2}
,"E"+(0,.5)*{3}
,"A"+(0,-.5)*{4}
,"B"+(0,-.5)*{5}
  \end{xy}
                            \qquad\strut&\qquad
\xymatrix{\bullet\ar@{-}[r]\POS!L(2)\drop{1}&\bullet\ar@{-}[r]\POS!D(1.5)\drop{2}
                                            &\bullet\ar@{-}[r]\POS!D(1.5)\drop{3}
                                              &\bullet\ar@{-}[r]\POS!D(1.5)\drop{4}
&\bullet\POS!R(1.5)\drop{5}
}
\end{array}
$$
\caption{Two examples of configurations represented by their nubs\\
  Left: $\I=\{\ve,\ 1,2,3,4,5,\ 13,15,23,24,25,34,45,\ 234\}$\\
  Right: $\I=\{\ve,\ 1,2,3,4,5,\ 13,14,15,24,25,35,\ 135\}$}
  \label{fig:qwpojsonmnmn}
\end{figure}

The configuration $\C$ is entirely determined by its collection of nubs. We shall use the convention of representing a configuration by depicting its collection of nubs, each nub being depicted by surrounding its elements. In case the nub has only two elements, we depict it by an edge connecting its two elements. For instance, on Figure~\ref{fig:qwpojsonmnmn}, we have depicted two configurations by theirs nubs and listed or described in the caption the corresponding \iset s.

\begin{definition}
  Let $\C=(V,\I,\D)$ be a \config.

  Two \iset s\ $x$ and $y$ are said to be \emph{parallel}, which we denote by $x\parallel y$, whenever $x\cap y=\emptyset$ and $x\cup y\in\I$. In this case we write $xy$ instead of~$x\cup y$, which is consistent with our word notations previously introduced.

  For every \iset\ $x$, the \emph{\config\ relative to~$x$} is the triple $\C^{\parallel x}=(V^{\parallel x},\I^{\parallel x},\D^{\parallel x})$ with $V^{\parallel x}=\{a\in V\tq x\parallel a\}$, and:
  \begin{align*}
    \I^{\parallel x}&=\bigl\{y\in\I\cap\P(V^{\parallel x})\tq x\parallel y\bigr\}
                      &\D^{\parallel x}=\P(V^{\parallel x})\setminus\I^{\parallel x}
  \end{align*}
\end{definition}

\begin{figure}
$$
\begin{array}{c|c}
  \fbox{\xymatrix{\bullet\POS!L(2)\drop{1}&\bullet\POS!R(1.5)\drop{2}\\
  \bullet\POS!L(2)\drop{4}&\bullet\POS!R(1.5)\drop{3}
}}\qquad\strut&\qquad
    \fbox{\xymatrix{\bullet\POS!L(2)\drop{1}&\bullet\POS!R(1.5)\drop{2}\\
  &\bullet\POS!R(1.5)\drop{3}
}}
\end{array}
$$
  \caption{Left, a \config\ $\C$ with its only nub depicted. Right, the relative configuration $\C^{\parallel a}$, with $a=4$ a vertex, with its only nub depicted.}
  \label{fig:qpojqwkqama}
\end{figure}
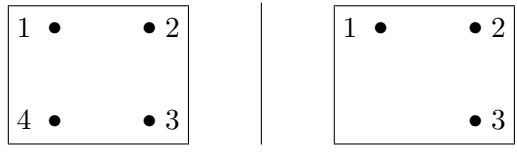

\begin{remark}
  \begin{enumerate}
  \item In general, $\I^{\parallel x}$ is a strict subset of $\I\cap\P(V^{\parallel x})$. Hence, $\I^{\parallel x}$ is not the mere restriction of $\I$ to~$V^{\parallel x}$. See an example below, illustrated on Figure~\ref{fig:qpojqwkqama}. The case where $\I^{\parallel x}=\I\cap\P(V^{\parallel x})$ corresponds to right-angled \config s, discussed below in~\S~\ref{sec:graphs-other-exampl}.
  \item We note that $\C^{\parallel\ve}=\C$.
  \item   Let $x\in\I$ and let $y\in\I^{\parallel x}$. Then $(\C^{\parallel x})^{\parallel y}=\C^{\parallel xy}$.
  \item Let $x,y\in\I$. Then it follows at once from the definitions of $x\leq y$ and of $z\parallel x$ that:
    \begin{gather}
      \label{eq:19}
      x\leq y\iff(\exists z\in\I^{\parallel x}\quad x=yz)
    \end{gather}
The \iset\ $z$ that appears in~\eqref{eq:19} is obviously unique, given by $z=y\setminus x$.
  \end{enumerate}
\end{remark}

\begin{example}
Let $V=[4]$ and $\I=\{x\in\P(V)\tq |x|\leq3\}$. The associated \config\ $\C=\SS_{4,3}$ is a particular case of the star configurations studied in more details in~\S\S~\ref{sec:star-config-s} and~\ref{sec:star-config-s-1}. The \config\  $\C$ is represented on the left of Figure~\ref{fig:qpojqwkqama}, where the only nub $V=[4]$ is depicted. For $a=4$, the relative \config\ $\C^{\parallel a}$ has the three vertices $1$, $2$,~$3$, and although $123$ is an \iset\ in~$\C$, it is not in~$\C^{\parallel a}$ since $1234\notin\I$. The relative \config\ $\C^{\parallel a}$ is depicted on the right of Figure~\ref{fig:qpojqwkqama}.
\end{example}

\subsection{Polynomials}
\label{sec:polynomials}

\begin{definition}[valuation]
  A \emph{valuation} on a \config\ $\C=(V,\I,\D)$ is a mapping $f:\I\to\RR_{>0}$ such that:
  \begin{gather}
    \label{eq:13}
f(\ve)=1\quad\text{and\quad}    (x,y\in\I)\quad x\parallel y\implies f(xy)=f(x)f(y)
\end{gather}

The \emph{uniform valuation} is the valuation given by $f(x)=1$ for all $x\in\I$.
\end{definition}

Obviously, the valuations of $\C$ are in bijection with~$(\RR_{>0})^V$; specifying $f$ on single letters, without any constraint other than being positive reals, is enough to determine $f$ entirely.

\begin{definition}
  A \emph{\wconfig} is a quadruple $(V,\I,\D,f)$ where $\C=(V,\I,\D)$ is a configuration, of which $f$ is a valuation. If the valuation is not specified, it is understood that the configuration is equipped with its uniform valuation.

  If $x\in\I$, we equip the relative \config\ $\C^{\parallel x}$ with the valuation~$f^{\parallel x}$, simply the restriction of $f$ to~$\I^{\parallel x}$.
\end{definition}

Let $\C=(V,\I,\D)$ be a configuration, and let $f:\I\to A$ be a mapping into an abelian group~$A$. It is well known that, for the partial order $\I$ which is a downward closed subset of $(\P(V),\subseteq)$, the Möbius transform~\cite{rota64} of~$f$ is the mapping $h:\I\to A$ defined by
\begin{gather}
  \label{eq:15}
  (x\in\I)\quad h(x)=\sum_{y\in\I\tq x\leq y}(-1)^{|y|-|x|}f(y)
\end{gather}
It allows to recover $f$ through the identity:
\begin{gather}
  \label{eq:17}
  (x\in\I)\quad f(x)=\sum_{y\in\I\tq x\leq y}h(y)
\end{gather}

Given a \wconfig\ $(V,\I,\D,f)$, we use the capital symbol $F$ to denote its extension $F:\I\to\RR[t]$ defined by $F(x)=t^{|x|}f(x)$ for $x\in\I$. Note that, evaluating the polynomials at a given positive real $t$ yields another valuation, namely $x\in\I\mapsto t^{|x|}f(x)$.

Let $H$ be the Möbius transform of~$F$, given by:
\begin{gather}
  \label{eq:18}
  (x\in\I)\quad H(x)=\sum_{y\in\I\tq x\leq y}(-1)^{|y|-|x|}t^{|y|}f(y)
\end{gather}
and, for $x\in\I$, let $H^{\parallel x}$ be the Möbius transform of~$F^{\parallel x}$. Since $f$ is a valuation on the one hand, and recalling the identity~\eqref{eq:19} on the other hand, one has:
\begin{gather}
  \label{eq:21}
  H(x)=F(x)H^{\parallel x}(\ve)
\end{gather}
The equation~\eqref{eq:21} is an identity between polynomials. This is a trivial identity for $x=\ve$.

\begin{definition}
  \label{def:4}
  The \emph{Möbius polynomial} of a \wconfig\ $(V,\I,\D,f)$ is the polynomial $H(\ve)$. An alternative notation for $H(\ve)$ is~$\mu(t)$, an alternative notation for $H^{\parallel x}(\ve)$ is $\mu^{\parallel x}(t)$.
\end{definition}

\begin{example}
  \begin{enumerate}
  \item The Möbius polynomial of the trivial \config\ $\emptyset$ is $\mu(t)=1$.
  \item For the \config s depicted on Figure~\ref{fig:qwpojsonmnmn}, the Möbius polynomials are respectively $\mu(t)=1-5t+7t^2-t^3$ (left) and $\mu(t)=1-5t+6t^2-t^3$ (right).
  \item For the \config\ depicted on Figure~\ref{fig:qpojqwkqama}, left, the various Möbius polynomials are: $H(\ve)=1-4t+6t^2-4t^3$, $H^{\parallel 4}(\ve)=1-3t+3t^2$, $H^{\parallel(43)}(\ve)=1-2t$, $H^{\parallel(432)}(\ve)=1$.
  \end{enumerate}
\end{example}

\section{Examples}
\label{sec:examples}

\subsection{Right-angled \config s}
\label{sec:graphs-other-exampl}

The configurations introduced above in Definition~\ref{def:2} can be interpreted as hypergraphs on their set of vertices. A special case is that of graphs, which we explore now through the notion of dependence relation. This notion appears in the literature dealing with trace monoids---also called free partially commutative monoids~\cite{cartier69}, heap monoids~\cite{viennot86} or right-angled Artin-Tits monoids~\cite{dehornoy14}.

A dependence relation $D$ on a set of vertices $V$ is defined as a binary, symmetric and reflexive relation on~$V$. The pair $(V,D)$ is called a \emph{dependence pair}. The associated independence relation $I$ is the binary, symmetric and irreflexive relation on $V$ defined by $I=(V\times V)\setminus D$. 

Given a dependence pair $(V,D)$, we define a configuration $\C=(V,\I,\D)$ by putting:
\begin{gather}
  \label{eq:22}
  \D=\bigl\{x\in\P(V)\tq \exists (a,b)\in V\times V\quad a\neq b\wedge \{a,b\}\subseteq x\bigr\}
\end{gather}
Equivalently, the set $\I$ is given by the set of cliques of the relation~$I$, ie:
\begin{gather}
  \label{eq:27}
x\in\I\iff\bigl(\forall (a,b)\in x\times x\quad a\neq b\implies (a,b)\in I\bigr) 
\end{gather}

\begin{definition}
  \label{def:3}
  A \config\ defined as in~\eqref{eq:22} is called a \emph{right-angled \config}.
\end{definition}

Examples include the complete graphs on a set of vertices~$V$, where all subsets are \dset s\ except for singletons and the empty set; and dimer models where $V=\{1,\ldots,n\}$ with $D=\{(i,i+1)\tq 1\leq i<n\}$, of which the example for $n=5$ was depicted on the right of Figure~\ref{fig:qwpojsonmnmn}.

If the \dset s in $\D$ are to be interpreted as global conflicts, then the following proposition shows that right-angled \config s are those configurations for which ``the conflicts are binary''; in other words, those for which the conflicts are localized at the scale of vertices.

\begin{proposition}
  \label{prop:1}
Let $\C=(V,\I,\D)$ be a \config. Then the following are equivalent:
\begin{enumerate}[\normalfont(i)]
  \item\label{item:1} $\C$ is a right-angled \config
  \item\label{item:2} All nubs of $\C$ have cardinal~$2$
  \item\label{item:3} For every \iset s  $x_1,\ldots,x_k\in\I$ such that $i\neq j\implies x_i\parallel x_j$, the set $x=x_1\cup\dots\cup x_k$ is an \iset
  \item\label{item:4} For every vertices $a_1,\ldots a_k$ such that $i\neq j\implies a_i\parallel a_j$, the set $x=\{a_1,\dots,a_k\}$ is an \iset
  \item\label{item:5} For every $x\in\I$, it holds that $\I^{\parallel x}=\I\cap\P(V^{\parallel x})$
  \item\label{item:6} For every vertex $a\in V$, it holds that $\I^{\parallel a}=\I\cap\P(V^{\parallel a})$
  \end{enumerate}

  Furthermore, if $\C$ is a right-angled \config, then all the relative \config s $\C^{\parallel x}$ for $x\in\I$ are also right-angled.
\end{proposition}

\begin{proof}
  \ref{item:1}$\implies$\ref{item:2}.\quad  Let $\C$ be a right-angled \config, induced by a dependence relation~$D$, and let $x$ be one of its nubs. Then $x\supset \{a,b\}$ for some pair $(a,b)\in D$ with $a\neq b$, and thus $x=\{a,b\}$ by minimality of~$x$.

  \ref{item:2}$\implies$\ref{item:3}.\quad Let $x_1,\ldots,x_k\in\I$ be such that  $i\neq j\implies x_i\parallel x_j$, and let $x=x_1\cup\dots\cup x_k$. Seeking a contradiction, assume that $x\notin\I$. Then $x\in\D$. Let $d\subseteq x$ be a nub. Then $|d|=2$ by assumption, say $d=\{a,b\}$. Then $a\in x_i$ and $b\in d_j$ for some distinct $i$ and~$j$, which contradicts that $x_i\parallel x_j$.

The implication \ref{item:3}$\implies$\ref{item:4} is trivial, \ref{item:4}$\implies$\ref{item:5} is immediate and \ref{item:5}$\implies$\ref{item:6} is trivial.

\ref{item:6}$\implies$\ref{item:1}. Let $d$ be a nub, let us prove that $|d|=2$. Seeking a contradiction, assume that $|d|>2$. For every two distinct vertices $a,b\in d$, one has that $\{a,b\}$ is a strict subset of~$d$. By minimality of~$d$, it implies that $\{a,b\}\in \I$, hence $a\parallel b$. By the assumption, it implies that the union of all these vertices is an \iset, hence $d\in\I$, which is a contradiction. Now set  $\widetilde D=\{(a,b)\in V\times V\tq \{a,b\}\in\D\}$ and let $D$ be the reflexive completion of~$\widetilde D$. Then $\C$ is the right-angled \config\ induced by~$D$.
\end{proof}


Let $\C$ be a right-angled \config, induced by the dependence pair $(V,D)$, and let $I$ be the associated independence relation. The Möbius polynomial polynomial $H(\ve)$ as defined in Definition~\ref{def:4} coincides with the classical Möbius polynomial $\mu(t)$ of the graph~$(V,I)$: \begin{gather}
  \label{eq:20}
  \mu(t)=\sum_{x\in\C}(-1)^{|x|}t^{|x|}
\end{gather}
where $\C$ denotes the set of cliques of the graph $(V,I)$, including the empty clique and the singletons sets. The polynomial $\theta(t)=\mu(-t)$ is the independence polynomial of the graph~$(V,D)$.

\subsection{Star \config s}
\label{sec:star-config-s}

\begin{definition}
  \label{def:5}
  Let $n,k\geq1$ be integers with $k\leq n$. The \emph{$(n,k)$-star \config} is the configuration $\SS_{n,k}=(V,\I,\D)$ where $V=[n]$, of size~$n$, and $\I=\{x\in\P(V)\tq |x|\leq k\}$. We also set $\SS_{0,0}=\emptyset$.
\end{definition}

$\SS_{n,k}$ is a right-angled \config\ if and only if $n=k=0$ or if $k=1$ or if $k=n$.

The Möbius polynomial of $\SS_{n,2}$ is $\mu_{n,2}(t)=1-nt+\frac{n(n-1)}2t^2$. More generally, the Möbius polynomial of $\SS_{n,k}$ is given by the polynomial $(1-t)^n$, truncated to the terms of degree at most~$k$. 

Let $a=n$, a vertex of $\SS_{n,k}$ with $k>1$. Then $(\SS_{n,k})^{\parallel a}=\SS_{n-1,k-1}$. More generally, for $x\in\I$ with $|x|<k$, the relative \config\ $(\SS_{n,k})^{\parallel x}$ is isomorphic to~$\SS_{n-|x|,k-|x|}$; whereas, if $|x|=k$ then $\SS_{n,k}=\emptyset$. Hence the character ``not right-angled'' is passed on to the relative \config s of~$\SS_{n,k}$, excepted for the trivial case when $x$ is of size~$k$ (see an illustration on Figure~\ref{fig:qpojqwkqama} for~$\SS_{4,3}$).

\section{Irreducibility}
\label{sec:irreducibility}

Let $V_1$ and $V_2$ be two subsets of $V$ such that $V_1\cap V_2=\emptyset$ and $V_1\cup V_2=V$. Let $\I_1=\{x\in\I\tq x\subseteq V_1\}$ and $\I_2=\{x\in\I\tq x\subseteq V_2\}$. There is obviously a one-to-one mapping $\varphi:\I\to\I_1\times\I_2$ defined by $\varphi(x)=(x\cap V_1,x\cap V_2)$. We defined the pair $(V_1,V_2)$ to form a \emph{decomposition} of $V$ if:
\begin{gather}
  \label{eq:16}
  \forall x\in\P(V)\quad x\in\I\iff(x\cap V_1\in\I)\wedge( x\cap V_2\in\I)
\end{gather}

This is equivalent to saying that $\varphi$ is onto, and thus bijective. The decomposition $(V_1,V_2)$ is said to be \emph{trivial} if $V_1=\emptyset$ or $V_2=\emptyset$.

\begin{definition}
  \label{def:1}
  A \config\  is \emph{reducible} if there exists a non trivial decomposition; it is
 \emph{irreducible} otherwise.
\end{definition}

\begin{remark}
  The fact that $\varphi$ is bijective implies at once the following product form for the Möbius polynomial $\mu(t)$ of~$V$: $\mu(t)=\mu_1(t)\mu_2(t)$, where $\mu_1$ and $\mu_2$ are respectively the Möbius polynomials of~$(V_1,\I_1)$ and of $(V_2,\I_2)$.
\end{remark}

Recalling the notion of nub from Definition~\ref{def:6}, we introduce the notion of nub path as follows.

\begin{definition}
  Given two vertices $a,b\in V$, a \emph{nub path from $a$ to~$b$} is a finite and non empty sequence $x_1,\ldots,x_k$ of nubs such that $a\in x_1$, $b\in x_k$ and $x_i\cap x_{i+1}\neq\emptyset$ for all $i=1,\dots,k-1$. The \emph{connected component of $a$} is the set of vertices $b$ that can be reached from $a$ through a nub path.
\end{definition}

\begin{lemma}
  \label{lem:1}
  Let $a_0\in V$, let $V_1$ be the connected component of~$a$ and let $V_2=V\setminus V_1$. Then $(V_1,V_2)$ is a decomposition of~$V$.
\end{lemma}

\begin{proof}
Let  $x_1\in\I_1$ and $x_2\in\I_2$, and let $x=x_1\cup x_2$. We prove that $x\in\I$; $x_1$~and $x_2$ being arbitrary, this will prove the claim of the lemma. Seeking a contradiction, assume that $x\notin\I$. Then $x$ contains a nub~$d$. Neither $d\subseteq x_1$ nor $x\subseteq x_2$ can hold since $x_1\in\I$ and $x_2\in\I$. Hence $d$ contains some vertex $a\in x_1$ and some vertex $b\in x_2$. But this implies that $b$ belongs to the connected component of~$a$, which is~$V_1$: this is a contradiction.
\end{proof}

\begin{proposition}
  A \config\  is irreducible if and only if it is connected.
\end{proposition}

\begin{proof}
  If a \config\  is irreducible, then Lemma~\ref{lem:1} implies at once that it is connected.

  Conversely, let $V$ be a connected \config. Seeking a contradiction, assume that it is reducible, and let $(V_1,V_2)$ be a non trivial decomposition of~$V$. Since $V$ is connected, there is a nub path that goes from $V_1$ to~$V_2$, and therefore there exists a nub $d$ satisfying $d\cap V_1\neq\emptyset$ and $d\cap V_2\neq\emptyset$. Then $d\cap V_1$ being a strict subset of~$d$, which is a nub, it must belong to~$\I$, and similarly $d\cap V_2\in\I$ for the symmetric argument. And since $(V_1,V_2)$ is a decomposition of~$V$, it implies that $d\in\I$, which is a contradiction.
\end{proof}

\begin{remark}
  If $\C$ is a right-angled \config\ induced by the dependence pair $(V,D)$, then its nubs are the edges of this graph. Hence $\C$ is irreducible if and only if $(V,D)$ is connected in the usual sense.
\end{remark}

\begin{example}
  The $(n,k)$-star \config\ is irreducible if and only if $n=k=1$ or if $k<n$.
\end{example}

\section{A formula for the derivative of the Möbius polynomial}
\label{sec:form-deriv-mobi}

\subsection{The general formula}
\label{sec:general-formula}

For a polynomial $P\in\RR[t]$, let $\ddt P$ denote the derivative of~$P$.

\begin{theorem}
  \label{thr:1}
  For every configuration $\C=(V,\I,\D)$, the derivative of its Möbius polynomial~$H(\ve)$ satisfies:
  \begin{gather}
    \label{eq:9}
    \ddt H(\ve)=-\sum_{a\in V}H^{\parallel a}(\ve)
  \end{gather}
\end{theorem}

\begin{proof}
  We proceed in two steps.

  \emph{First step}. We observe that there is a unique family of integers $\alpha_{\gamma,\theta}$ defined for $\gamma,\theta\in\I$ such that $\gamma\parallel\theta$, which satisfies the following identities:
    \begin{gather}
      \label{eq:6}
     (\gamma\parallel\theta)\qquad \alpha_{\gamma,\theta}=|\theta|-\Biggl(\sum_{\substack{(\delta,\lambda)\in\I\times\I\tq           \delta\neq\ve\\
\delta\parallel\gamma,\ \lambda\parallel\ (\gamma\delta),\ \delta\lambda=\theta
}}\alpha_{\gamma\delta,\lambda}\Biggr)
    \end{gather}

    Indeed, if it exists, it must satisfy:
    \begin{gather}
      \label{eq:8}
(\gamma\in\I)\qquad  \alpha_{\gamma,\ve}=0
\end{gather}
and the terms are then defined by a descending induction with respect to~$\lambda$, starting from the maximal~$\lambda$; for those, $\alpha_{\gamma,\ve}=0$ are the only valid terms by the maximality of~$\lambda$.

 Next, we observe that $\alpha_{\gamma,\theta}$ are actually given by:
 \begin{gather}
\label{eq:7}
   \alpha_{\gamma,\theta}=
   \begin{cases}
     1,&\text{if $|\theta|=1$}\\
     0,&\text{otherwise}
   \end{cases}
 \end{gather}
 The proof of~\eqref{eq:7} is by increasing induction on~$|\theta|$. The case $|\theta|=0$ follows from~\eqref{eq:8}. The case $|\theta|=1$ derives from~\eqref{eq:6}; for, with $|\theta|=1$, the sum in the right hand side of~\eqref{eq:6} contains only one term, which is~$0$. If $|\theta|=2$, say $\theta=ab$ with $(a,b)\in V\times V$, then:
 \begin{gather*}
   \alpha_{\gamma,ab}=2-\alpha_{\gamma a,b}-\alpha_{\gamma b,a}=0
 \end{gather*}
by the preceding case. We now proceed to the induction case for $|\theta|=k>2$, assuming that $\alpha_{\gamma,\lambda}=0$ if $1<|\lambda|<k$. By this induction hypothesis, and since $\alpha_{\gamma\delta,a}=1$ if $a\parallel(\gamma\delta)$ by the previous case, $\alpha_{\gamma,\theta}$~reduces to:
\begin{gather*}
  \alpha_{\gamma,\theta}=|\theta|-\#\bigl\{
(\delta,a)\in\I\times V\tq\delta\parallel\gamma,\ a\parallel(\gamma\delta),\ \delta a=\theta\}
  \bigr\}=0
\end{gather*}
completing the induction step and the proof of~\eqref{eq:8}. This completes our first step.

\medskip
\emph{Second step.} We now aim to prove the following:
\begin{gather}
  \label{eq:10}
  (\gamma\in\I)\qquad \ddt H^{\parallel \gamma}(\ve)=-\sum_{\lambda\in\I\tq\lambda\neq\ve,\ \lambda\parallel\gamma}\alpha_{\gamma,\lambda}t^{|\lambda|-1}
H^{\parallel \gamma\lambda}(\ve)
\end{gather}
Specializing~\eqref{eq:10} for $\gamma=\ve$, and combined with~\eqref{eq:7}, this will prove~\eqref{eq:9} and complete the proof of the theorem.

We prove~\eqref{eq:10} by descending induction on~$\gamma$. If $\gamma$ is maximal in~$\I$, then $H^{\parallel \gamma}(\ve)=1$, hence $\ddt H^{\parallel \gamma}(\ve)=0$, which coincides with the right-hand side of~\eqref{eq:10} since the sum is empty when $\gamma$ is maximal.

To prove the induction step, it is enough to prove~\eqref{eq:10} for $\gamma=\ve$ assuming that it holds for all $\gamma\neq\ve$. Hence we are brought to prove:
\begin{gather}
  \label{eq:11}
  \ddt H(\ve)=-\sum_{\lambda\in\I\tq\lambda\neq\ve}\alpha_{\ve,\lambda}t^{|\lambda|-1}H_\lambda(\ve)
\end{gather}

Deriving the identity $\sum_{\gamma\in\I}H(\gamma)=1$, and using the identity $H(\gamma)=t^{|\gamma|}H^{\parallel \gamma}(\ve)$ for $\gamma\neq\ve$, we get:
\begin{gather}
  \label{eq:12}
  \ddt H(\ve)=-\sum_{\gamma\in\I,\ \gamma\neq\ve}\Bigl(
|\gamma|t^{|\gamma|-1}H^{\parallel \gamma}(\ve)+t^{|\gamma|}\ddt H^{\parallel \gamma}(\ve)
  \Bigr)
\end{gather}

The induction assumption provides allows to substitute the right-hand side of~\eqref{eq:10} in~\eqref{eq:12}, yielding:
\begin{align}
\notag
  \ddt H(\ve)&=-\sum_{\gamma\in\I,\ \gamma\neq\ve}\Bigl(
  |\gamma|t^{|\gamma|-1}H^{\parallel \gamma}(\ve)\Bigr)+
  \sum_{\substack{(\lambda,\gamma)\in\I\times\I\tq\\ \gamma\neq\ve,\ \lambda\neq\ve,\
  \lambda\parallel\gamma}}\Bigl(
t^{|\gamma\lambda|-1}\alpha_{\gamma,\lambda} H^{\parallel \gamma\lambda}(\ve)
  \Bigr)\\
  \label{eq:14}           &=-\sum_{\gamma\in\I, \gamma\neq\ve}
R_\gamma               t^{|\gamma|-1}H^{\parallel \gamma}(\ve)
\end{align}
where $R_\gamma$, obtained by regrouping the terms in the second sum of the right-hand side of~\eqref{eq:14}, is given by:
\begin{align*}
\notag  R_\gamma&=|\gamma|-\sum_{\substack{(\delta,\theta)\in\I\times\I\tq\\
  \delta\neq\ve,\ \theta\neq\ve,\ \delta\parallel\theta,\ \delta\theta=\lambda
  }}\alpha_{\delta,\theta}\\
\notag  &=|\gamma|-\sum_{\substack{(\delta,\theta)\in\I\times\I\tq\\
  \delta\neq\ve,\ \delta\parallel\theta,\ \delta\theta=\lambda
  }}\alpha_{\delta,\theta}&&\text{since $\alpha_{\delta,\ve}=0$}\\
  &=\alpha_{\ve,\lambda}&&\text{by~\eqref{eq:6}}
\end{align*}

We substitute $R_\gamma=\alpha_{\ve,\lambda}$ in~\eqref{eq:14} to obtain~\eqref{eq:11}, as expected. The proof of the theorem is complete.
\end{proof}

\subsection{A combinatorial application}
\label{sec:comb-appl}

Before we come to our targeted probabilistic applications, we mention a combinatorial application of the above formula for the derivative of the Möbius polynomial.

\begin{proposition}
  \label{prop:2}
  Let $\C$ be a \config\ with the following symmetry property:
  \begin{gather}
    \label{eq:23}
    (x,y\in\I)\quad |x|=|y|\implies\bigl(
\text{$\C^{\parallel x}$ and $\C^{\parallel y}$ are isomorphic}
    \bigr)
  \end{gather}
For $j\geq0$, let:
\begin{gather*}
(|x|=j)\quad \eta_j=\# V^{\parallel x}\,,  \qquad N_j=\#\{x\in\I\tq |x|=j\}
\end{gather*}
Then:
\begin{gather}
  \label{eq:24}
  (k\geq0)\quad N_k=\frac1{k!}\eta_0\dots\eta_{k-1}
\end{gather}
\end{proposition}

\begin{proof}
  Let $a_1\ldots a_K$ be an \iset\ of maximal size. Applying successively Theorem~\ref{thr:1}, and using the symmetry property~\eqref{eq:23}, one has:
  \begin{align*}
    \ddt H(\ve)&=-\eta_1 H^{\parallel a_1}(\ve)    &
                                                     \ddt H^{\parallel {a_1}}(\ve)&=-\eta_2 H^{\parallel a_1a_2}(\ve)\\
&\dots&\ddt H^{\parallel a_1\dots a_j}(\ve)&=-\eta_{j}H^{\parallel a_1\dots a_{j+1}}(\ve)
  \end{align*}
Evaluating the above identities at $t=0$ yields the successive derivatives of the polynomials at zero. Therefore:
\begin{gather}
  \label{eq:25}
  H^{\parallel a_1\dots a_{K-j}}(\ve)=1-\eta_{k-j}t+(\eta_{k-j}\eta_{k-j+1})\frac{t^2}{2}+\dots+(-1)^j(\eta_{K-j}\dots\eta_{K-1})\frac{t^j}{j!}
\end{gather}
and finally, for $j=K$ in the above:
\begin{gather}
  \label{eq:26}
  H(\ve)=1-\eta_0t+\eta_0\eta_1\frac{t^2}2+\dots+(-1)^{K}\eta_0\dots\eta_{K-1}\frac{t^K}{K!}
\end{gather}
But the coefficients of the polynomial $H(\ve)$ are also given, with a sign, by the number of \iset s of size~$j$, whence the identity~\eqref{eq:24}.
\end{proof}

For $V=[n]$ equipped with $\I=\P(V)$, the formula~\eqref{eq:24} is simply the formula for the binomial coefficients. 

In general, the identity~\eqref{eq:24} is trivial for $k=0,1$ and for $k>K$. It is already non trivial for $k=2$ as it says: $N_2=\frac12 N_1\eta_1$, since $\eta_0=N_1$.

Consider the set $V$ of vertices of a dodecahedron, so that $N_1=20$, and equip $V$ with a \config\ $(V,\I,\D)$ where $x\in\I$ if and only if $x$ is a subset of an edge of the dodecahedron (this is not a right-angled \config). Then $\eta_1=3$ since a given vertex belongs to three different edges. Hence the number of edges is $N_2=\frac12\times 20\times 3=30$, in case we had forgotten it.

\section{Probabilistic \config s}
\label{sec:prob-conf-s}

\begin{definition}
  \label{def:8}
  A \emph{probabilistic \config} is given by a quintuple $(V,\I,\D,f,t)$ where $(V,\I,\D,f)$ is a \wconfig\ and $t$ is a non negative real such that, for some probability space $(\Omega,\F,\pr)$, there exists a mapping $\varphi:V\to\F$ with the following properties:
  \begin{enumerate}
  \item For every $a\in V$: $\pr(\varphi(a))=tf(a)$
  \item For every $x\in\I$, the family $(\varphi(a))_{a\in x}$ is a family of independent events
  \item For every $x\in\D$, the family $(\varphi(a))_{a\in x}$ is a family of mutually exclusive events modulo~$0$, ie, $\pr\bigl(\bigcap_{a\in x}\varphi(a)\bigr)=0$.
  \end{enumerate}

  The probability space $(\Omega,\F,\pr)$ is said to be \emph{configured} by $(V,\I,\D)$. The \emph{rest} of the configured space is the probability
  \begin{gather}
    \label{eq:28}
    R=1-\pr\Bigl(\bigcup_{a\in V}\varphi(a)\Bigr)
  \end{gather}
\end{definition}

\begin{proposition}
  \label{prop:3}
  Let $\C=(V,\I,\D,f)$ be a \wconfig. A real $t\geq0$ makes $(\C,t)$ probabilistic if and only if:
  \begin{gather}
    \label{eq:29}
    \forall x\in\I\quad H(x)[t]\geq0
  \end{gather}

  In this case, a configured probability space is the finite set $X=\I$, together with the mapping $\varphi:a\to X$ defined by $\varphi(a)=\{x\in X\tq a\leq x\}$, and equipped with the unique probability distribution $m$ satisfying $m(\varphi(a))=tf(a)$ for all $a\in V$.

  For every other configured space $(\Omega,\F,\pr)$, there is a canonical random variable $(\Omega,\F)\to X$ with law~$m$; and the rest is necessarily given by $R=H(\ve)[t]$. 

In particular, the property that the family of events $(\varphi(a))_{a\in X}$ is a covering modulo $0$ of\/~$\Omega$ does not depend on~$\Omega$, but only on $\C$ and~$t$. This property is satisfied if and only if $t$ is a root of~$H(\ve)$.
\end{proposition}

\begin{proof}
  We are in the framework of our preliminary section, \S~\ref{sec:free-models-prob}. Recall the notations introduced therein, namely $\K$ for our purpose. Fix $t\geq0$. For every $x\in\K$, let $q_x=t^{|x|}f(x)$ if $x\in\I$ and let $q_x=0$ if $x\in\D$. Then $(\C,t)$ is probabilistic if and only if $(q_x)_{x\in\K}$ can be realized in the sense of~\S~\ref{sec:free-models-prob}.

  Hence, for every $x\in\K$, let $p_x=\sum_{y\in\K\tq x\leq y}(-1)^{|y|-|x|}q_y$. According to Proposition~\ref{MMprop:2}, the realization exists if and only if all $p_x$ are non negative. But $p_x=H(x)[t]$, and thus the condition~\eqref{eq:29} is necessary and sufficient for $(\C,t)$ to be probabilistic. The remaining of the statement is a reformulation of Proposition~\ref{MMprop:2}.
\end{proof}

\begin{theorem}
  \label{thr:2}
  Let $\C=(V,\I,\D,f)$ be a non trivial \wconfig. The set of non negative reals making $(\C,t)$ probabilistic is of the form $[0,t_0]$, where $t_0>0$ is given by:
  \begin{gather}
    \label{eq:30}
    t_0=\min\{t\geq0\tq\exists x\in\I \quad H(x)[t]=0\}
  \end{gather}

Furthermore, the rest $R(t)$ is decreasing on $[0,t_0]$, hence $t_0$ is also the real minimizing the rest $R(t)$, for $t$ ranging over the possible reals making $(\C,t)$ probabilistic.
\end{theorem}

\begin{proof}
Let:
\begin{align*}
  T&=\{t\geq0\tq \text{$(\C,t)$ is probabilistic}\}\\
  T_1&=\{t\geq0\tq \forall x\in\I\quad H(x)[t]\geq0\}
\end{align*}

Proposition~\ref{prop:3} shows that $T=T_1$.
  There exists at least a positive real $t$ and an \iset\ $x$ such that $H(x)(t)=0$. Indeed, let $x_1$ be an element of maximal size in~$\I$; $|x_1|>0$ since $\C$ is assumed to be non trivial. Pick a subset $x\leq x_1$  with $|x_1|=|x|-1$ and let $a_1$ be the vertex such that $x_1=xa_1$. Then $\{a,b\}\in \D^{\parallel x}$ for any $a,b\in V^{\parallel x}$ with $a\neq b$; indeed, otherwise one would have $xab\in\I$, contradicting the maximality of~$|x_1|$. Therefore, $H(x)=f(x)t^{|x|}H^{\parallel x}(\ve)=f(x)t^{|x|}(1-kt)$ where $k=\sum_{a\in V^{\parallel x}}f(a)$. Furthermore $k\geq f(a_1)>0$, hence $t=\frac1k$ is a positive root.

  Since $H(x)[0]=1$ for all $x\in\I$, the real $t_0$ defined by~\eqref{eq:30} exists and is positive. Furthermore, $H(x)[t]\geq0$ for all $x\in\I$ and for all $t\in[0,t_0]$. Henceforth:
  \begin{gather}
    \label{eq:32}
    [0,t_0]\subseteq T_1
  \end{gather}

Consider now a real $t>0$ such that $H(x)[t]\geq0$ for all $x\in\I$. Then we claim:
  \begin{itemize}
  \item[$(\dag) $] There exists $\delta>0$ such that $H(x)[r]\geq0$ for all $x\in\I$ and for all $r\in[t-\delta,t]$.
  \end{itemize}

To prove the claim~$(\dag)$, let $\delta\in(0,t)$ be such that all the polynomials $H(x)=f(x)t^{|x|}H^{\parallel x}(\ve)$, for $x$ ranging over~$\I$, are of constant sign and are either non decreasing or non increasing on~$[t-\delta,t]$. The same property then holds for all the polynomials $H^{\parallel x}(\ve)$. 

Seeking a contradiction, assume that $H(x)[t-\delta]<0$ for some $x\in\I$, and thus $H^{\parallel x}(\ve)[t-\delta]<0$, and pick an element $x$ maximal in $\I$ with this property. Since $H^{\parallel x}(\ve)[t]\geq0$ by hypothesis, and since $H^{\parallel x}(\ve)$ is monotonic and negative on~$[t-\delta,t)$, it must be increasing on~$[t-\delta,t]$. But, by Theorem~\ref{thr:1}, $\ddt H^{\parallel x}(\ve)=-\sum_{a\in V^{\parallel x}}H^{\parallel xa}(\ve)$, and all the terms in this sum are non negative on $[t-\delta,t]$ by the maximality property of~$x$. Hence $\ddt H^{\parallel x}(\ve)$ is non increasing, which is a contradiction and completes the proof of the claim~$(\dag)$.
  
According to standard analysis arguments, it follows from~$(\dag)$ just proved that $T_1=[0,t_1]$ for some positive real~$t_1$. From~\eqref{eq:32}, we conclude that $t_0\leq t_1$.

  We now prove the following claim:
  \begin{itemize}
  \item[$(\ddag)$] There exists $x\in\I$ and a real $\delta>0$ such that  $H^{\parallel x}(\ve)[r]<0$ for all $r\in(t_0,t_0+\delta)$.
  \end{itemize}

  Indeed, let $x\in\I$ be maximal among those $x\in\I$ such that $H^{\parallel x}(\ve)[t_0]=0$. Then $H^{\parallel y}(\ve)[t_0]>0$ for all $y\in\I$ with $x\neq y$ and $x\leq y$. By Theorem~\ref{thr:1}, $\ddt H^{\parallel x}(\ve)=-\sum_{a\in V^{\parallel x}}H^{\parallel xa}(\ve)[t_0]>0$, therefore $H^{\parallel x}(\ve)$ is strictly decreasing on a neighborhood of~$t_0$, which completes the proof of the claim~$(\ddag)$.

  The statement~$(\ddag)$ implies that $t_0<t_1$ is impossible, hence $t_0=t_1$.

  It remains only to prove that $t\mapsto R(t)$ is decreasing on $[0,t_0]$. By Theorem~\ref{thr:1}, $\ddt R(t)=-\sum_{a\in V}H^{\parallel a}(\ve)[t]$ and all the terms are non negative on $[0,t_0]$ since $T=T_1$ on the one hand, and since $H(a)[t]=f(a)tH^{\parallel a}(\ve)[t]$ on the other hand. The proof is complete.
\end{proof}

\begin{definition}
  \label{def:7}
  Let $\C$ be a non trivial \wconfig\ and let $t_0$ be the real defined by~\eqref{eq:30}. We say that $t_0$ is the \emph{critical root} of~$\C$. Let $R=H(\ve)[t_0]$ be the rest of the associated probabilistic \config. We say that $\C$ is:
  \begin{enumerate}
  \item of type \one\ if $R=0$
  \item of type \two\ if $R>0$.
  \end{enumerate}
\end{definition}

\begin{remark}
  In other words, a non trivial configuration is of type \one\ if and only it can be realized as a configured space which is covered modulo~$0$ by the elementary events determined by the vertices.
\end{remark}

The following result is a direct consequence of Theorem~\ref{thr:2}.

\begin{theorem}
  \label{thr:3}
  A non trivial \wconfig\ is of type~\one\ if and only if there exists a real $t>0$ such that $H(\ve)[t]=0$ and $H(x)[t]\geq0$ for all $x\in\I$. If it exists, this real is unique and is the critical root of the \wconfig.
\end{theorem}

\section{Examples, continued}
\label{sec:examples-continued}

\subsection{Right-angled \config s}
\label{sec:right-angled-config}

As we already mentioned in the Introduction, the study of right-angled \config s greatly benefits of a combinatorial interpretation of their Möbius polynomials using associated trace monoids. Let $\C=(V,\I,\D,f)$ be a weighted right-angled \config\ induced by a dependence pair $(V,D)$ and let $I=(V\times V)\setminus D$. The trace monoid associated with $\C$ is then $\M=\langle V\;|\; ab=ba\text{ for all $(a,b)\in I$}\rangle$. The generating series of $(\M,f)$ is
\begin{gather}
  \label{eq:31}
  G(t)=\sum_{x\in\M}f(x)t^{|x|}
\end{gather}

A result of~\cite{cartier69,viennot86} is the formal identity $G(t)\mu(t)=1$, where $\mu(t)=\sum_{x\in\I}(-1)^{|x|}f(x)t^{|x|}$ is the Möbius polynomial of the weighted monoid, or of the weighted configuration. This identity becomes valid in the field of reals for all $t\in(0,r)$, where $r<+\infty$ is the radius of convergence of~$G(t)$.

Hence the series $G(t)$ is rational, inverse of~$\mu(t)$. Being rational and with non negative coefficients, its radius of convergence $r$ is one of its poles, and thus a root of~$\mu(t)$, and $\mu(z)$ has no other root in the complex plane of smallest modulus than~$r$. It actually holds that $\mu(z)$ has no other root of modulus $r$ than $r$ itself. Furthermore, $\mu(t)G(t)=1$ implies that $\mu(t)>0$ for all $t\in[0,r)$.

Now for every $x\in\I$, the relative \config\ $C^{\parallel x}$ is also a right-angled \config, and the associated trace monoid $\M^{\parallel x}$ is nothing but a sub-monoid of~$\M$, namely the submonoid $\M^{\parallel x}$ generated by all vertices $a\parallel x$. As a consequence, its generating series $G^{\parallel x}(t)$ satisfies evidently $G^{\parallel x}(t)\leq G(t)$ for every $t\in[0,r^{\parallel x})$ where $r^{\parallel x}$ is the radius of convergence of~$G^{\parallel x}(t)$. Henceforth, $r^{\parallel x}\geq r$ on the one hand, and:
\begin{gather}
  \label{eq:33}
( x,y\in\I,\quad t\in(0,r])\quad x\leq y\implies \mu^{\parallel x}(t)\leq\mu^{\parallel y}(t)
\end{gather}

As a consequence, for $t=r$, one has simultaneously $\mu(t)=0$ and $\mu^{\parallel x}(t)\geq0$. By Theorem~\ref{thr:3}, it follows that $\C$ is of type~\one, and that $r$ is its critical root.

A thorough study of trace monoids, using in particular the normal form for their elements introduced in~\cite{cartier69}, and an appeal to Perron-Frobenius theory for primitive matrices, allows to prove some additional results in $\C$ is irreducible. Firstly, the root $r$ is a simple root of~$\mu(t)$~. An secondly, still assuming that $\C$ is irreducible, then it actually holds that $r<r^{\parallel x}$ for all $x\in\I$, $x\neq\ve$. Therefore, for the critical root $t=r$, we have $\mu^{\parallel x}(t)>0$ for all $x\neq\ve$. 

We have obtained the following result (details about the facts stated above can be found in~\cite{krob03,goldwurm00,csikvari13,viennot86,cartier69}).

\begin{proposition}
  \label{prop:4}
  Every right-angled \config\ $\C$ is of type~\one. Its critical root is the only root of smallest modulus in the complex plane of its Möbius polynomial. If $\C$ is irreducible, this root $r$ is simple and satisfies $\mu^{\parallel x}(r)>0$ for every $x\in\I$, $x\neq\ve$.
\end{proposition}

\subsection{Star \config s}
\label{sec:star-config-s-1}

It is elementary to establish the following result (proof omitted).

\begin{proposition}
  \label{prop:5}
  For $n,k\geq1$ with $k\leq n$, the types of the star configurations $\SS_{n,k}$ follow the following pattern:\
  \begin{center}
    \normalfont
    \begin{tabular}{cc|ccccccc}
      &&1&2&3&4&5&6&$\xrightarrow{\ k\ }$\\
      \hline
      &1&\one\\
      &2&\one&\one\\
      &3&\one&\two&\one\\
      &4&\one&\two&\one&\one\\
      &5&\one&\two&\one&\two&\one\\
$n\downarrow$      &6&\one&\two&\one&\two&\one&\one
    \end{tabular}
  \end{center}
  
  In particular, $\SS_{n,n-1}$ is of type \one\ if $n$ is even, and of type \two\ if $n$ is odd. The critical root for $\SS_{n,n-1}$ is always~$\frac12$.
\end{proposition}

\begin{figure}
\begin{gather*}
  \begin{xy}
   0;/r1.7cm/:
     (0,0)="G",
    *\xycircle(1,1){}
,"G";/r1.1cm/:
"G",
{\xypolygon4"A"{~><{@{}}~={45}{}}},
"A1"*{ab},"A2"*{ac\phantom{b}},"A3"*{bc},
"G";/r2.2cm/:
"G",
{\xypolygon4"B"{~><{@{}}~={45}{\frac14}}},
"G",
{\xypolygon4"C"{~><{@{}}~={0}{}}},
"C1";"C3"**@{-},
"C2";"C4"**@{-}
\end{xy}
\hspace{.15\textwidth}
  \begin{xy}
    0;/r1.7cm/:
     (0,0)="G",
     *\xycircle(1,1){}
     ,"G";/r1.1cm/:
     "G",
{\xypolygon8"A"{~><{@{}}~={22.5}{}}},
     "A1"*{a},"A2"*{abd},"A3"*{abc},"A4"*{acd},"A5"*{bcd},"A6"*{b},"A7"*{c},"A8"*{d},
/r2.2cm/:
"G",
{\xypolygon8"B"{~><{@{}}~={22.5}{\frac18}}}
,"G",{\xypolygon8"C"{~><{@{}}~={0}{}}},
"C1";"C5"**@{-},
"C2";"C6"**@{-},
"C3";"C7"**@{-},
"C4";"C8"**@{-}
\end{xy}
\end{gather*}
\caption{The configured disk for $\SS_{n,n-1}$: left, with $n=3$, of type \two\ and with minimal rest~$\frac14$; right, with $n=4$ and of type \one.}
  \label{fig:ppopoiqwqwqw}
\end{figure}
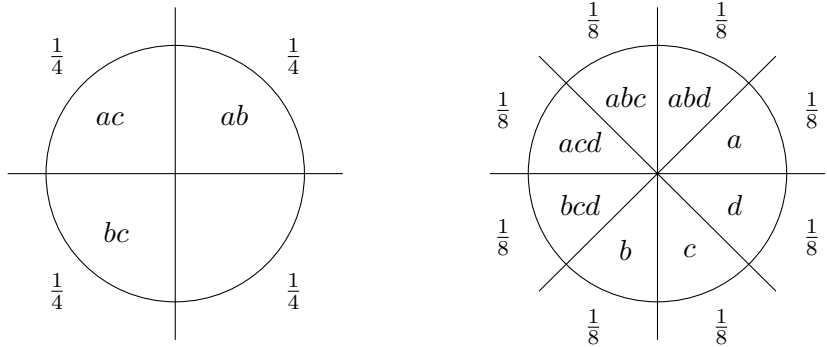

We depict on Figure~\ref{fig:ppopoiqwqwqw} an illustration of the configured space for $\SS_{n,n-1}$ with $n=3,4$, taking for $\Omega$ a disk of area $1$ and at the critical root~$t=\frac12$. Observe on each picture that all the various independence and exclusion conditions are met.

\section{Further questions}
\label{sec:further-questions}

We point out two questions that seem to be of interest, for the challenge that they represent and for their potential applications.

The first question concerns the classification of irreducible \config s. Is there some criterion to determine whether a given irreducible \wconfig\ is of type \one\ or of type~\two~? And more specifically, is it true that the type of an irreducible \wconfig\ is independent of the valuation~$f$, and is thus purely a property of the underlying \config? Of course, the computability of such a criterion, if it exists, should be examined as well.

Another question concerns the comparison between configured spaces, relatively to a given \wconfig, and other probability spaces. A very nice result of Shearer's paper~\cite{shearer85} is the following. Let $(V,I)$ be an independence graph, let $\C$ be the associated \config, say of critical root~$t_0$, and let $R(t)$ be the rest of the associated canonical configured space, for $t\in[0,t_0]$. Assume now be given a probability space $(\Omega,\F,\pr)$ with events corresponding to the vertices of~$V$, and such that these events satisfy the independence structure prescribed by $(V,I)$. Hence, $(\Omega,\F,\pr)$ is almost a realization of~$\C$, except that it is not required that families of events in $\D$ are mutually exclusive in~$\Omega$. Let also $\alpha(t)=1-\pr\bigl(\bigcup_{a\in V}a\bigr)$. Then it is shown that $\alpha(t)\geq R(t)$.

Obtaining the same inequality for general configurations is hopeless. Indeed, for $\C=\SS_{3,2}$ for instance, with minimal rest~$\frac14$, it is very much possible to construct a probability space entirely covered by three events and satisfying the pairwise independence conditions: take simply $\Omega=\{*\}$ and $a=b=c=\{*\}$. This is precisely because one requires that $\{a,b,c\}$ is mutually exclusive that the rest must be positive.

Hence, obtaining comparison results based on general configured spaces requires, if this is achievable, to set up some additional conditions that remain to be found.

\paragraph*{Graphics.}
\label{sec:graphics}

The pictures of this paper have been typeset with ~\Xy-pic, a graphic package for \TeX\ designed by K.~Rose and R.~Moore.

\paragraph*{Remerciements.}
\label{sec:aknowldegements}

Je voudrais remercier mes collègues et amis A.~Benveniste, J.~Mairesse, V.~Jugé, S.~Gouëzel et A.~Durand. Ils ont chacun su m'encourager et m'épauler à un moment ou à un autre de ma vie mathématique.

\printbibliography

\end{document}